\documentclass[11pt]{article}
\usepackage{amssymb,amsmath,amsthm,hyperref,verbatim,fullpage}

\newtheorem{theorem}{Theorem}

\newtheorem{conjecture}[theorem]{Conjecture}

\newcommand{\F}{\mathcal{F}}
\newcommand{\G}{\mathcal{G}}

\newcommand{\Int}{\mathrm{Int}}
\newcommand{\ucstar}{\textrm{UC}^*}
\begin{document}
\title{Union-closed families with small average overlap densities.}
\author{David Ellis\footnote{School of Mathematics, University of Bristol.}}
\date{3rd December 2020}
\maketitle

\begin{abstract}
\noindent In this very short paper, we show that the average overlap density of a union-closed family $\mathcal{F}$ of subsets of $\{1,2,\ldots,n\}$ may be as small as $\Theta((\log_2 \log_2 |\F|)/(\log_2 |\F|))$, for infinitely many positive integers $n$.
\end{abstract}

\section{Introduction}

If $X$ is a set, a family $\mathcal{F}$ of subsets of $X$ is said to be {\em union-closed} if the union of any two sets in $\F$ is also in $\F$. The celebrated Union-Closed Conjecture (a conjecture of Frankl \cite{duffus}) states that if $X$ is a finite set and $\F$ is a union-closed family of subsets of $X$ (with $\F \neq \{\emptyset\}$), then there exists an element $x \in X$ such that $x$ is contained in at least half of the sets in $\F$. Despite the efforts of many researchers over the last forty-five years, and a recent Polymath project \cite{polymath} aimed at resolving it, this conjecture remains wide open. It has only been proved under very strong constraints on the ground-set $X$ or the family $\F$; for example, Balla, Bollob\'as and Eccles \cite{bbe} proved it in the case where $|\F| \geq \tfrac{2}{3} 2^{|X|}$; more recently, Karpas \cite{karpas} proved it in the case where $|\F| \geq (\tfrac{1}{2}-c)2^{|X|}$ for a small absolute constant $c >0$; and it is also known to hold whenever $|X| \leq 12$ or $|\F| \leq 50$, from work of Vu\v{c}kovi\'c and \v{Z}ivkovi\'c \cite{vz} and of Roberts and Simpson \cite{rs}.

In 2016, a Polymath project \cite{polymath} was convened to tackle the Union-Closed Conjecture. While it did not result in a proof of the conjecture, several interesting related conjectures were posed. Among them was the `average overlap density conjecture'.

If $X$ is a finite set and $\F \subset \mathcal{P}(X)$ with $\F \neq \emptyset$, we define the {\em abundance} of $x$ (with respect to $\F$) by $\gamma_x = |\{A \in \F:\ x \in A\}|/|\F|$, i.e., $\gamma_x$ is the probability that a uniformly random member of $\F$ contains $x$. A natural first quantity to consider, in trying to prove the Union-Closed Conjecture, is the average abundance of a uniformly random element of the ground set, i.e., $\mathbb{E}_{x \in X}[\gamma_x]$; if this quantity were always at least $1/2$, the Union-Closed Conjecture would immediately follow. A moment's thought shows that this is false, however, e.g.\ by considering the union-closed family $\{\emptyset,\{1\},\{1,2,3\}\} \subset \mathcal{P}(\{1,2,3\})$, which has average abundance $4/9$. Similarly, for any $n \in \mathbb{N}$, the union-closed family $\F = \{\emptyset,\{1\},\{1,2\},\ldots,\{1,2,\ldots,\lfloor \sqrt{n}\rfloor\},\{1,2,3,\ldots,n\}\} \subset \mathcal{P}(\{1,2,\ldots,n\})$ has average abundance $\Theta(1/\sqrt{n}) = \Theta(1/|\F|)$.

It is natural to consider the expected abundance of a random element of the ground-set $X$ chosen according to other (non-uniform) distributions on $X$. The following was considered in the Polymath project \cite{polymath}. We define the {\em average overlap density} $\text{AOD}(\F)$ of $\F$ to be the expected value of $\gamma_x$, where $x$ is a uniformly random element of a uniformly random nonempty member of $\F$:
\begin{align}
\label{eq:reform}
\text{AOD}(\F) &:= \frac{1}{|\F \setminus \{\emptyset\}|} \sum_{A \in \F \setminus \{\emptyset\}} \frac{1}{|A|} \sum_{x \in A}\gamma_x\nonumber \\ 
& = \frac{1}{|\F \setminus \{\emptyset\}|} \sum_{A \in \F \setminus \{\emptyset\}} \frac{1}{|A|} \sum_{x \in A}\frac{|\{B \in \F:\ x \in B\}|}{|\F|} \nonumber\\
& = \frac{1}{|\F \setminus \{\emptyset\}|} \sum_{A \in \F \setminus \{\emptyset\}}\left( \frac{1}{|\F|} \sum_{B \in \F} \frac{|A \cap B|}{|A|}\right) \nonumber\\
& = \mathbb{E}_{A \in \F \setminus \{\emptyset\}} \mathbb{E}_{B \in \F}\left[\frac{|A \cap B|}{|A|}\right].
\end{align}
(The first and second expectations in (\ref{eq:reform}) are of course over a uniformly random element of $\F \setminus \{\emptyset\}$, and a uniformly random element of $\F$, respectively.) The last equality justifies the `average overlap' terminology. The average overlap density conjecture stated that if $X$ is a finite set, and $\F$ is a union-closed family of subsets of $X$ with $\F \neq \emptyset$ and $\F \neq \{\emptyset\}$, then the average overlap density of $\F$ is at least $1/2$. Clearly, it would immediately imply the Union-Closed Conjecture. 

Unfortunately, the average overlap density conjecture was quickly shown to be false (during the Polymath project \cite{tb}); an infinite sequence of union-closed families $\F_n \subset \mathcal{P}(\{1,2,\ldots,n\})$ was constructed with $\text{AOD}(\F_n) = 7/15+o(1)$ as $n \to \infty$. However, the following weakening of the average overlap density conjecture remained open.
\begin{conjecture}
\label{conj:adocw}
There exists an absolute positive constant $c>0$ such that the following holds. Let $n \in \mathbb{N}$ and let $\F \subset \mathcal{P}(\{1,2,\ldots.n\})$ be union-closed with $\F \neq \{\emptyset\}$. Then the average overlap density of $\F$ is at least $c$.
\end{conjecture}
Conjecture \ref{conj:adocw} would immediately imply the weakening of the Union-Closed Conjecture where $1/2$ is replaced by the absolute positive constant $c$.

In this very short paper, we prove the following.
\begin{theorem}
\label{thm:main}
For infinitely many positive integers $n$, there exists a union-closed family $\F$ of subsets of $\{1,2,\ldots,n\}$ whose average overlap density is $\Theta((\log_2 \log_2 |\F|)/(\log_2 |\F|))$.
\end{theorem}
This disproves Conjecture \ref{conj:adocw} in a strong sense. It follows from an old result of Knill \cite{knill} that if $\F \subset \mathcal{P}(\{1,2,\ldots,n\})$ is union-closed, then there exists $x \in \{1,2,\ldots,n\}$ with abundance $\gamma_x = \Omega(1/(\log_2 |\F|))$, so the average overlap density can, in the best-case scenario, only be used to improve this lower bound by a factor of $\Theta(\log_2 \log_2 |\F|)$.

\section{Proof of Theorem \ref{thm:main}}
For $n \in \mathbb{N}$, we write $[n]: = \{1,2,\ldots,n\}$ for the standard $n$-element set, and if $\mathcal{G} \subset \mathcal{P}(X)$, the {\em union-closed family generated by $\G$} is defined to be the smallest union-closed family of subsets of $X$ that contains $\G$.

Let $k,m,s \in \mathbb{N}$ with $s \leq k-2$ and $m \geq 2$, and let $n=km$. Partition $[n]$ into $m$ sets $B_1,\ldots,B_m$ with $|B_i| = k$ for all $i$; in what follows, we will refer to the $B_i$ as `blocks'. For each $i \in [m]$, choose a subset $T_i \subset B_i$ with $|T_i|=s$, and let $T = \cup_{i=1}^{m}T_i$. Now let $\F \subset \mathcal{P}([n])$ be the union-closed family generated by $\{B_i \cup \{j\} :\ i \in [m],\ j \in T\}$. Note that every set in $\F$ contains at least one block. The number of sets in $\F$ containing exactly one block is $m2^{(m-1)s}$, and in general, for each $j \in [m]$, number $N_j$ of sets in $\F$ containing exactly $j$ blocks is ${m \choose j} 2^{(m-j)s}$, so
$$N: = |\F| = \sum_{j=1}^{m} N_j = 2^{(m-1)s} \sum_{j=1}^{m} {m \choose j} 2^{-(j-1)s}.$$
For each $j \in [m]$, define $p_j: = N_j/N$; this is of course the probability that a uniformly random member of $\F$ contains exactly $j$ blocks. We note that 
$$\frac{p_{j+1}}{p_{j}} = \frac{N_{j+1}}{N_j} = \frac{m-j}{j+1}2^{-s} \leq m2^{-s} \quad \forall j \in [m-1].$$
 Write $\tau : = m2^{-s}$. For any $x \in [n] \setminus T$, we clearly have
 $$\gamma_x = \frac{1}{m}\sum_{j=1}^{m} jp_j,$$
since the conditional probability that $x$ is contained in a random member $A$ of $\F$, given that $A$ contains exactly $j$ blocks, is $j/m$. We have $p_j \leq \tau^{j-1}p_1$ for all $j \in [m]$, and therefore for any $x \in [n] \setminus T$, we have
$$\frac{1}{m} \leq \gamma_x \leq \frac{1}{m}(1+2\tau + 3\tau^2 + \ldots + m\tau^{m-1}) \leq \frac{1}{m}(1+4\tau) \leq \frac{2}{m},$$
provided $\tau = m2^{-s} \leq 1/4$. Now, every member $A$ of $\F$ contains at least one block, so for any member $A$ of $\F$, the probability a uniformly random element of $A$ is in $T$, is at most $\frac{ms}{k}$. Crudely, we have $1/2 \leq \gamma_x \leq 1$ for all $x \in T$, since $A \mapsto A \cup \{x\}$ is an injection from $\{A \in \F:\ x \notin A\}$ to $\{A \in \F:\ x \in A\}$, for any $x \in T$. Hence, we have
\begin{equation}
\label{eq:bound1}
\frac{1}{m} \leq \text{AOD}(\F) \leq \left(1-\frac{ms}{k}\right)\cdot  \frac{2}{m} + \frac{ms}{k} \cdot 1 \leq \frac{2}{m} + \frac{m^2 s}{n},
\end{equation}
again provided $\tau = m2^{-s} \leq 1/4$. Now we wish to minimize the right-hand side of (\ref{eq:bound1}), subject to the constraint $m2^{-s} \leq 1/4$; clearly the optimal choice is to take $s = \lceil \log_2 m \rceil + 2$, which yields
\begin{equation}
\label{eq:bound2}
\frac{1}{m} \leq \text{AOD}(\F) \leq \frac{2}{m} + \frac{m^2 \log_2 m}{n} + O(m^2/n).
\end{equation}
It is clear that the optimal choice of $m$ to minimize the right-hand side of (\ref{eq:bound2}) is
$$m = \Theta\left(\left(\frac{n}{\log_2 n}\right)^{1/3}\right),$$
yielding $\text{AOD}(\F) = \Theta(((\log_2 n)/n)^{1/3})$. Since, with these choices, we have
$$\log_2|\F| = \Theta(n^{1/3} (\log_2 n)^{2/3}),$$
it follows that
$$\text{AOD}(\F) = \Theta\left(\frac{\log_2 \log_2 |\F|}{\log_2 |\F|}\right),$$
proving Theorem \ref{thm:main}.

We proceed to note two further properties of the above construction. Firstly, the average abundance of a uniformly random element of $[n]$ (with respect to $\F$) satisfies
$$\mathbb{E}_{x \in [n]} [\gamma_x] = \Theta\left(\frac{\log_2 \log_2 |\F|}{\log_2 |\F|}\right).$$

Secondly, the family $\F$ constructed above does not separate the points of $[n]$. (We say a family $\F \subset \mathcal{P}([n])$ {\em separates the points of $[n]$} if for any $i \neq j \in [n]$ there exists $A \in \F$ such that $|A \cap \{i,j\}|=1$. It is easy to see that, in attempting to prove the Union-Closed Conjecture, we may assume that the union-closed family in question separates the points of the ground set, and this assumption was adopted for much of the Polymath project \cite{polymath}.) However, it is easy to see that the union-closed family $\F \cup \{[n] \setminus \{j\}:\ j \in [n]\}$ has asymptotically the same average overlap density as $\F$ (and asymptotically the same average abundance as $\F$), and does separate the points of $[n]$.

\end{document}